\documentclass[11pt]{article}

\usepackage{amssymb,amsmath}

\def\diag{\mathop{{\rm diag}}\nolimits}
\def\eq{\Leftrightarrow}
\def\hs{\hspace*{\parindent}}
\def\id{\mathord{{\rm Id}}}

\def\int{\mathop{{\rm Int}}\nolimits}

\def\m{\medskip}
\def\proof{{\hs \bf Proof.\ }}
\def\qed{\hfill $\Box$}

\def\sl{{\rm SL}}
\def\stab{\mathop{{\rm Stab}}\nolimits}

\def\vs{\vspace{\baselineskip}}

\def\<{\mathopen{<}}
\def\>{\mathclose{>}}

\def\R{\mathord{\mathbb R}}

\def\apl#1#2#3#4{\begin{eqnarray*}#1 & \to & #2\\ #3 & \mapsto & #4 \end{eqnarray*}}

\def\dtto#1{\big\downarrow\rlap{$\scriptstyle #1$}}
\def\tto#1{\mathop{\longrightarrow}\limits^{{#1}}}

\def\eq{\mathrel{\Leftrightarrow}}

\newtheorem{theo}{\bf \hs Theorem}[section]
\newtheorem{defn}[theo]{\bf \hs Definition}
\newtheorem{prop}[theo]{\bf \hs Proposition}
\newtheorem{lemma}[theo]{\bf \hs Lemma}
\newtheorem{corol}[theo]{\bf \hs Corollary}

\def\a{\mathord{\mathfrak a}}
\def\Ad{\mathop{\rm Ad}}
\def\c{\mathord{\mathfrak c}}
\def\cont{\mathop{\rm Cont}}
\def\diag{\mathop{\rm diag}}
\def\g{\mathord{\mathfrak g}}

\def\hs{\hspace*{\parindent}}
\def\id{\mathord{\rm Id}}
\def\k{\mathord{\mathfrak k}}
\def\p{\mathord{\mathfrak p}}
\def\proof{\hspace*{\parindent}{\bf Proof. }}
\def\qed{\hfill $\Box$}

\def\sl{\mathord{\rm SL}}
\def\so{\mathord{\rm SO}}
\def\stab{\mathop{\rm Stab}\nolimits}
\def\Tr{\mathop{{\rm Tr}}}
\def\uin{\mathrel{\underline\in}}
\def\vs{\medskip}
\def\w{\mathord{\mathfrak w}}

\def\tto#1{\mathop{\longrightarrow}\limits^{{#1}}}

\begin{document}
\title{A Topological Approach to Unifying Compactifications of Symmetric Spaces}
 
\author
{Pedro J.\ Freitas \\
Center for Linear and Combinatorial Structures\\
Department of Mathematics --- Faculty of Sciences\\
University of Lisbon\\
{\tt pedro@ptmat.fc.ul.pt}}

\date{June 23, 2011}
\maketitle

\begin{abstract}
In this paper we present a topological way of building a compactification 
of a symmetric space from a compactification of a Weyl Chamber.

\end{abstract} 

\section{Introduction}

\hs There has been some recent interest in finding ways to unify the processes of obtaining compactifications of symmetric spaces $G/K$ (see \cite{GJT}, \cite{BJ}), where $G$ is a semisimple connected non-compact Lie group with finite center, and 
$K$ is a maximal compact subgroup. These unifying procedures use mainly concepts from differential geometry or Lie group theory, and aim at producing general ways to obtain known compactifications of symmetric spaces, such as the Visual, maximal Satake, maximal Furstenberg, Martin and Karpelevi\v c compactifications. 

In \cite{GJT}, it is shown that these compactifications actually depend on the compactification of a flat through a given point $o=K\in G/K$, and on the fact that they are $K$-equivariant. These properties, along with another property on the compactification of intersection of Weyl chambers, actually identify the compactification. Some of these constructions have a shortcoming, they do not allow for a natural $G$-action, a problem that was overcome in \cite{BJ}, with a different approach to these general constructions, this time making more use of topology, along with parabolic groups.

In this paper, we also present a topological way of building a compactification of $X=G/K$ from a compactification of the Weyl chamber centered at $o$, generalizing the constructions that were done \cite{GJT} with flats, for each of the known  compactifications, listed above. We prove some properties about this compactification, including existence and uniqueness, in a rather general setting, requiring only an extra condition on the compactification of the Weyl chamber. We then identify some known compactifications as particular cases of this construction. We have the same shortcoming of not being able to define a $G$-action, but the setting in which we work is quite general

As an addendum to this, we present a different way of building a compactification of a symmetric space, using generalized Busemann functions. We establish that it is indeed a compactification and make some conjectures on how to obtain the known compactifications in this manner. 

We would like to thank Vadim Kaimanovich for presenting the problem, and for many suggestions and discussions on this subject. \medskip

\section{General Concepts}

We start by defining the notation (either well known or taken from \cite{GJT}, with minor adjustments) and the concepts necessary.  The results that follow can be found in \cite{GJT} and \cite{He}.

We recall that we take $G$ to be a semisimple connected non-compact Lie group with finite center, and 
$K$ be a maximal compact subgroup. Denote by $\g$ and $\k$ the Lie algebras of $G$ and $K$
respectively.

Let $\g=\k\oplus \p$ be the Cartan decomposition of $\g$, 
$\p$ being the orthogonal complement of $\k$ in $\g$, with respect to the Killing form $B$. 
The space $\p$ can be identified with the tangent space to $X$ at the coset $K$, which we'll
denote by $o$. The restriction of the Killing form to this space is positive definite, and
thus provides an inner product in $\p$. 

We take $\a$ to be a fixed Cartan subalgebra of $\p$, $\a^+$ a fixed Weyl chamber,
$\Sigma$ the set of all the roots of $\g$ with respect to $\a$ (the so-called restricted roots),
$\Sigma^+$ the set of positive roots, $\Delta$ the set of the simple roots. We denote by $d$ be the rank of $G$ (the dimension of $\a$).
\vs

The action of $G$ on $X$ is by left multiplication. The adjoint action of
$K$ on $\p$ is given by the derivative of $X\mapsto g\exp(tX)g^{-1}$ at $t=0$, for $g\in K$ and
$X\in K$. These actions are in correspondence to each other, meaning that
$\exp((\Ad g)(X))=g \exp(X)$. 

Every element of $G$ can be written as $k.\exp(X)$ with $k\in K$ and $X\in \p$---this is an easy consequence of the Cartan decomposition, which states that every element of $G$ can be expressed as $k_1\exp(X)k_2$, $k_1,k_2\in K$, $X\in \overline{\a^+}$.\medskip

We write $A^+:=\exp(\a^+)$, $\w:=\overline{\a^+}$ and $W:=\exp(\w)=\exp (\overline{\a^+})$. From this we can say that every point in $X$ can be presented as $p=kx.o$, $k\in K$, $x\in W$, which means that the $K$-orbit of $W$ is whole symmetric space. Moreover, the element $x\in W$ is uniquely defined, and is called the generalized radius of $p=kx.o$ with respect to $o$. The element $k$ is unique modulo the stabilizer of $x$ for the action of $K$ over $W$.
\vs

As an example, take $G=\sl(d+1,\R)$. In
this particular case, $K=\so (d+1, \R)$, the Killing form on $\g$, the set of all
matrices of trace zero, is given by $B(M,N)=2(d+1)\Tr(MN^t)$. The set $\k$ is the set
of all skew-symmetric matrices, and $\p$ is the set of all symmetric matrices of
trace zero, and thus $B$ restricted to $\p$ is just a multiple of the usual
scalar product. The Cartan subalgebra $\a$ is the set of all diagonal matrices with
positive entries placed in strictly increasing order. The set of simple roots is:
$$\Delta=\{L_1-L_2,L_2-L_3,\ldots, L_d-L_{d+1}\},$$
where for each $i$, $L_i$ is the form dual to the matrix $E_{ii}$, the matrix with
all entries equal to zero except the $(i,i)$ entry, which is equal to one. The set of all roots
is $\Sigma=\{L_i-L_j,\ 1\leq i\neq j \leq d+1 \}$, and the set of all positive roots is 
$\Sigma^+=\{L_i-L_j, 1\leq i < j \leq d+1\}$.\medskip

Given a topological group $H$, we'll say that a topological space $T$ is an {\em
$H$-space} if there is an action of $H$ on $T$ (which we'll denote by a dot) and the
map 
\apl{H\times T}{T}{(h,a)}{h.a}
is continuous. If $R$ is another $H$-space, and $\phi:T\to R$ is a continuous map,
we say that $\phi$ is {\em $H$-equivariant} if, for any $h\in H, a\in T$,
$\phi(h.a)=h.\phi(a)$.

If $R$ is compact, $\phi$ is an embedding, and $\phi(T)$ is dense in $R$, we'll say
that
$(\phi,R)$ (or simply $R$ if there is no confusion about the map involved) is a {\em
compactification} of $T$. If $\phi$ is $H$-equivariant, we'll say that $R$ is an
{\em $H$-compactification}.

\section{Building a Compactification from the Weyl Chamber}

We are now concerned with the definition of a compactification of the space $X$
via compactifications of the closed Weyl chamber. There are a few compactifications of
$X$ that can be presented this way, such as the compactifications of Furstenberg, Satake, Karpelevi\v{c} and Martin.\medskip



Now, we will build for a $K$-invariant compactification of $X$ that once restricted to $W$ will be $\tilde{W}$. \medskip

{\bf Weyl chamber faces}. Let $I\subseteq \Delta$. Adjusting the definition and properties in \cite[p.\ 25]{GJT}, we define a {\em Weyl chamber face} as
$$\begin{array}{rcl}
\c_I & = & \{H\in \overline{\a^+} : \alpha(H)>0 \text{ if and only if }Ê\alpha \notin I\}\\
& = & \{H\in \overline{\a^+} : \alpha(H)=0 \text{ if and only if }Ê\alpha \in I\}
\end{array}
$$
and $C_I:=\exp(\c_I)$ (in \cite{GJT} the sets $C_I$ were contained in $\overline{\a^+}$). 

The Weyl chamber faces constitute a partition of the closed Weyl chamber, since they are pairwise disjoint and their union is the closed Weyl Chamber---we note, for instance, that $\exp(\a^+)=C_\emptyset$ and $o=C_\Delta$.\medskip

Given a face $\c_I$ of the Weyl chamber $W$, denote by $C_K(\c_I)$ the centralizer of $\c_I$ in $K$: $k\in C_K(C_I)$ if and only if $k\in K$ and for all $x\in \c_I$, $\Ad_k(x)=x$. We denote by $\bar \c_I$ the closure of $\c_I$ in $\a$. We have corresponding definitions for $C_I$.\medskip

We denote by $\tilde C_I$ the compactification we get for $C_I$, restricted from $\tilde W$.\medskip

We now intruduce an extra rquirement for the compactification $\tilde W$. It is known that if $kx.o = ry.o$, for $k,r\in K$ and $x,y\in W$, then we must have $x=y$ --- see \cite[Th. 1.1, p. 420]{He}. By the same theorem, if $x\in \exp(\a^+)$, then $k^{-1}r$ has to be in the center of $G$. If, however, $x=\exp(H)$, with $H\in \overline{a^+} \setminus \exp(\a^+)$, then it must lie in a Weyl chamber face $\c_I$, and we must have that $\Ad_{k^{-1}r} H=H$. By Lemma 3.10 and Proposition 2.15 in \cite{GJT}, if the element $k^{-1}r$ fixes an element in $\c_I$, it must centralize (that is, pointwise fix) the whole face $\c_I$. Therefore, for $x,y\in W$, we can say that 
$kx.o = ry.o$ if and only if $x=y$ and $k^{-1}r$ fixes the Weyl chamber face $C_I$ such that $x\in C_I$.\medskip

Now we want the compactification $\tilde W$ to satisfy a similar property. However, one cannot expect an element of $\tilde W\setminus W$ to belong to the compactification of only one Weyl chamber face. Thinking strictly about closure in $\a$, it is easy to check, by looking at the definition of $\c_I$, that  
$$\overline{\c_I} = \bigcup_{J\supseteq I} \c_J,$$
and similarly for $C_I$. Therefore, for each $x\in \overline{A^+}$, the set
$$\{J: x\in \overline{C_J} \}$$
has a maximum, which is exactly the set $I$ such that $x\in C_I$. 

This will thus be of the properties we will demand of the compactification $\tilde W$.\medskip

From now on, we will assume that 
\begin{itemize}
\item $\tilde W$ is metrizable and 
\item for each $x\in \tilde W$, the set $\{ J: x\in \tilde C_J\}$ has a maximum. If $I$ is this maximum, we will write $x \uin \tilde C_I$. 
\end{itemize}

We note that for $x\in \overline{A^+}$, $x\in C_I \eq x\uin \tilde C_I$. We will say that a compactification of $W$ satisfying these conditions is {\em facially stratified}.\medskip 

\medskip

{\bf The equivalence relation}. For $I\subset \Delta$, denote by $\stab(I)$ the centralizer of $C_I$ in $K$, which coincides with the centralizer of any point in $C_I$, as we have seen (again, see Lemma 3.10 and Proposition 2.15 in \cite{GJT} for a description of this set). \medskip

We note that, just by checking definitions, we have
$$I\subseteq J \eq C_I \supseteq C_J \eq \stab(I) \subseteq \stab(J).$$

We also note that $\stab(I)$ is a closed set (to see this, one can use the definition of $\stab(I)$ or the description in Proposition 2.15 of \cite{GJT}). Since it is a closed subset of a compact set, it must be compact.\medskip 

Now consider the space $K\times W$, and the map $\pi_1:K\times W\to X$ defined naturally by
$\pi_1(k,x):=kx.o$. Consider now the compact space $K\times \tilde{W}$ and
its quotient by the relation $\sim$ defined by the following rule: 
\begin{quote}
For $k,r\in K$, $x,y \in \tilde{W}$, $(k,x)\sim (r,y)$ if and only if $x=y$ and
if $x\uin C_I$ then $k^{-1}r\in \stab(I)$.
\end{quote}It is easy to check that it is an equivalence relation, under the conditions we have for the compactification $\tilde W$. 

Moreover, from what we have seen, for $x,y\in W$, $(k,x)\sim (r,y)$ if and only if $kx.o=ry.o$. If $kx.o=ry.o$, then we must have $x=y$ and if $x\in C_I$, then $x\uin \tilde C_I$ and we must have $k^{-1}r C_I =C_I$, which means $k^{-1}r\in \stab(I)$. The converse is equally simple. 

This which allows us to identify the set $(K\times W)/\sim$ with $X$. Therefore, this equivalence relation states that generalized radii must exist in $(K\times \tilde W)/\sim$.\medskip

Denote by ${\tilde X}$ the quotient space endowed with the quotient topology. Now take the inclusion and projection maps
$$\iota_1:K\times W \to K\times \tilde{W} \qquad \pi_2:K\times \tilde{W}\to
\tilde X.$$ 
By what we have said, the following diagram commutes.
$$\begin{array}{ccc}
K\times W & \tto{\pi_1} & X\\
\dtto{\displaystyle\iota_1} &  & \dtto{\displaystyle\iota}\\
K\times \tilde{W} & \tto{\pi_2} & {\tilde X}
\end{array}$$

It is clear that $\iota$ is the identity map onto $\iota(X)$, so once we prove that $\iota(X)$ is dense in $\tilde X$ and that $\tilde X$ is compact, we will have that $\tilde X$ is a compactification of $X$. We start by proving that $\tilde X$ is metrizable. 

\begin{prop} The map $\pi_2$ is closed.
\label{pi2_closed}
\end{prop}

\proof By theorem 10, p.\ 97 in \cite{Ke}, this is equivalent to showing that if a set $M\subset K\times \tilde W$ is closed, then 
$$\sim [M] :=\{z\in K\times \tilde W : z\sim z' \text{ for some } z'\in M \}$$
is closed. 

Take $M\subseteq K\times \tilde W$, a closed set. Since $K\times \tilde W$ is compact, the $M$ is also compact and hence both projections of $K$ and on $\tilde W$ must also be compact. 

To prove that $\sim[M]$ is closed, take a converging sequence $(k_n,x_n)\in \mbox{$\sim[M]$}$, with $(x_n,k_n)\to (k,x)$. We wish to show that $(k,x)\in \sim [M]$. 

We must have $(k_n,x_n)\sim (r_n,x_n)$, with $(r_n,x_n)\in A$ (the first coordinate has to be equal, according to the definition of $\sim$). Since the first projection of A is compact, we can take a converging subsequence of $r_n$ --- we'll consider that $r_n$ is already convergent to $r$, to simplify notation. Since $A$ is closed, we must have $(r,x)\in A$. 

Let $I\subseteq \Delta$ be such that $x\uin \tilde C_I$, so that $x\in \tilde C_J \Rightarrow J\subseteq I$. Then we must have that, for $n$ large enough, $x_n\in \tilde C_J$ for some $J\subseteq I$, and $k_n^{-1}r_n \in \stab (J) \subseteq \stab (I)$. Since $\stab (I)$ is closed, we must have $k^{-1}r\in \stab(I)$. Therefore $(k,x)\sim (r,x)$, with $(r,x)\in A$, so $(k,x)\in \sim[M]$, as we wished.\qed

\begin{theo} The space $\tilde X$ is metrizable and compact. It is a compactification of $X$.\label{compact_metrizable}
\end{theo}
\proof By the corollary of Theorem 20, p. 148 and Theorem 12 of p.\ 99 of \cite{Ke}, if $\pi_2$ is closed and the classes for $\sim$ are compact, then $\tilde X$ is metrizable. From what we have seen, the class of $(k,x)$, for $x\uin C_I$ is $k\stab(I)\times \{x\}$, which is clearly a compact set. Since we just proved that $\pi_2$ is closed, we have metrizability. 

The space $\tilde X$ is clearly compact, since $K\times \tilde W$ is compact and $\pi_2$ is continuous. To see that $\iota(X)$ is dense in $\tilde X$, take $(k,x)/\sim \in \tilde X$, we have that there is a sequence $(k_n,x_n)\in K\times W$ converging to $(k,x)$, and by continuity of $\pi_2$, we must also have convergence in $\tilde X$, which finishes the proof.\qed\medskip

{\bf A $K$-action}. It is now easy to see that we have a continuous action of $K$ on $\tilde X$, naturally defined as $r.(k,x)/\sim := (rk,x)/\sim$. It is well defined, since if $(k_1,x)\sim (k_2,x)$, then, if $x\uin \tilde C_I$, $k_1^{-1}k_2\in \stab(I)$ and
$$(rk_1)^{-1}(rk_2)=k_1^{-1}r^{-1}rk_2=k_1^{-1}k_2\in \stab(I)$$
and $(rk_1,x)\sim (rk_2,x)$. 

In view of this, from now on, for $(k,x)/\sim \in \tilde X$, we will denote $(k,x)/\sim$ by $kx.o$. We finish this section with three important properties of the compactification $\tilde X$.

\begin{prop} \label{props_tildew} The compactification $\tilde X$ has the following properties:\smallskip

1. It is a $K$-compactification.\smallskip

2. The compactification of $W$ considered as a subset of $\tilde X$ is $\tilde W$. \smallskip

3. The compactification $\tilde X$ respects intersections of Weyl chambers, that is, for $k,r\in K$, 
$$k\tilde W \cap r\tilde W = \widetilde{k W \cap rW}.$$
\end{prop}
\proof 1. To see that the $K$-action is continuous, and that $\tilde X$ is a $K$-space, consider the following diagram. 

$$\begin{array}{ccc}
K\times (K\times \tilde W) & \tto{\kappa'} & K\times \tilde W\\
\dtto{\id \times \pi_2} &  & \dtto{\pi_2}\\
K\times \tilde{X} & \tto{\kappa} & {\tilde X}
\end{array}$$

We denoted by $\kappa$ the action of $K$ on $\tilde X$ that we have just defined, and by $\kappa'$ the map $(r,(k,x))\to (rk,x)$. We wish to see that $\kappa$ is continuous, which, according to Theorem 9, p. 95 of \cite{Ke}, is equivalent to saying that $\kappa \circ (\id\times \pi_2)$ is continuous. Since
$$\kappa \circ (\id\times \pi_2) = \pi_2\circ \kappa'$$
we have the desired continuity and $\tilde X$ becomes a $K$-space. 

2. The image $\pi_2(\{\id\}\times \tilde W)$ is clearly homeomorphic to $\tilde W$ and is the compactification of $W$ considered as a subset of $\tilde X$.

3. Since $kW\cap rW \subseteq k\tilde W \cap r\tilde W$, and the second set is closed in $\tilde X$, we must have 
$\widetilde{k W \cap rW} \subseteq k\tilde W \cap r\tilde W$.

Conversely, let $kx.o=rx.o\in k\tilde W \cap r\tilde W$. If $x\uin \tilde C_I$, we have that $k^{-1}r\in  \stab(I)$, so $kC_I=rC_I\subseteq kW\cap rW$. Since $kx.o\in k\tilde C_I$, $kx.o\in \widetilde{kW\cup rW}$, as we wished.\qed\medskip

{\bf Examples}. There are a few known compactifications that are particular cases of our
compactification $\tilde X$, originating from different compactifications of $W$, namely, the visual, maximal Furstenberg, maximal Satake, Karpelevi\v c and Martin compactifications. We refer to descriptions given in \cite{GJT} and prove that the compactification of the Weyl chamber is, in each case, facially stratified.\medskip

 For the visual compactification, restricted to the Weyl chamber, we can associate each point in the boundary with a unit vector $v\in W$ (see p.\ 23). If $v\in C_I$, then $v\uin C_I$ by the structure of the faces of the Weyl chamber, and $v$ is fixed by $\stab I$. \medskip
 
The dual cell compactification, which is isomorphic to the maximal Satake compactification (Theorem 4.43) and the maximal Furstenberg compactification (Theorem 4.53) is described in page 41, definition 3.35. We have, in the notation used in this definition, that $(C_I(\infty),a^I)\uin \tilde C_I$ if $a^I=0$. If $a^I\neq 0$, then $x\uin C_\emptyset=A^+$. In any case, the limit point is fixed by $\stab(I)$. \medskip

The formal limits for the Karpelevi\v c compactification of $W$ are described in Definition 5.14, and the action of $K$ on these limits is given on p.\ 85. As in the previous case, if $H^I=0$, then the set $I$ appearing in the definition of the formal limit $x$ determines the Weyl chamber wall $C_I$ for which $x\uin C_I$, if $H^I\neq 0$, then $x\uin C_\emptyset=A^+$. Again, this limit is preserved by $\stab(I)$. \medskip

 For the (most general) Martin compactification, the limits are described in Theorem 8.2 and Proposition 8.20. According to this last proposition, the points $x_{I,a,L}\in \tilde W$ depend of three parameters: $I\subseteq \Delta$, $a\in C_I^\perp$, $L\in \overline{C_I}$ with $||L||=1$. Turning to the discussion about $I$-directional sequences on Proposition 8.9, it is easy to conclude that $x_{I,a,L}\in \tilde C_I$ if and only if $J\subseteq I$ and $a=0$. Hence $x_{I,0,L}\uin \tilde C_I$, and again according to Proposition 8.20, this limit point is preserved by $\stab(I)$. 

We note that $I$-directional sequences (defined on p.\ 119), which are used here, are the $C_I$-fundamental sequences (defined on p.\ 35), which are the ones used for the dual cell compactification, with a limiting direction $L$. This reflects the fact that the Martin compactification is a refinement of the dual cell compactification.

\section{Uniqueness}

We now recall the concept of fundamental subsequence, taken from \cite{GJT}. 

\begin{defn}
Let $X$ be a non-compact topological space, and $\bar X$ a compactification. A set of
sequences ${\cal C}$ of $X$ is called a {\em system of fundamental sequences} (for $\bar X$) if
\begin{itemize}
\item all sequences in ${\cal C}$ are convergent in $\bar{X}$, and
\item every sequence in $X$ has a subsequence in ${\cal C}$.
\end{itemize}
\end{defn}

{\bf Example}. For any $K$-equivariant compactification of $X$, then we can take
as a set of fundamental sequences, the set 
$$\{k_nx_n:k_n {\rm\ and\ } x_n{\rm\ converge} \}.$$
This is very easy to verify. To start with, these sequences have to converge
because the action of $K$ is continuous. Now, given any sequence
$r_ny_n$, there is a converging subsequence of $r_n$, say $r_{\alpha_n}$,
because $K$ is compact, and then there is a converging subsequence of $x_{\alpha_n}$
in the restriction of the compactification to $W$. Thus we find a fundamental
subsequence of any sequence in $X$.\medskip

{\bf Remark}. {\em In a $K$-compactification, not all convergent sequences are necessarily
fundamental}. Take, for instance, the one-point compactification of ${\mathbb
H}^2$---consider the upper half plane model. Then any sequence $k_n(2n).i$ converges
to infinity, no matter which sequence $k_n\in \so(2)$ we choose. 

For a more refined example, take the visual compactification of the symmtric space $\sl(3,\R)/\so(3)$,
with the point $o=\so(3)$. Take the sequence $x_n:=\diag(-n,-n,2n)\in W$, and
$k_n:=\diag((-1)^n,(-1)^n,1)$. Then $k_n\exp(x_n).o=\diag(e^{-n},e^{-n},e^{2n}).o$, with
converging limit direction given by the vector $\diag(0,0,1)\in W$. Thus, the
sequence converges in the visual compactification. Notice, however, that the
sublimits of $k_n$ are in the stabilizer of the limit direction.\medskip

Even though not all convergent sequences are fundamental, still, fundamental
sequences, along with their respective limits, determine the sequences of $X$ which converge in $\bar X$.

\begin{prop}
Let there be given a set of fundamental sequences for a compactification $\bar X$ of a
space $X$. Then a sequence $x_n$ in $X$ converges to $x\in\bar X$ if and only if every
fundamental subsequence of $x_n$ converges to $x$. 
\label{FundConv}
\end{prop}
\proof If $x_n$ converges to $x$, then obviously, every subsequence converges to $x$. 
Conversely, assume that every fundamental subsequence of $x_n$ converges to $x$,
and suppose that $x_n$ doesn't converge to $x$. Then, there must exist a
neighborhood of $x$, $U$, and subsequence of $x_n$ that remains outside $U$. Taking
now a fundamental subsequence of this subsequence, we have that, under our
assumption, it must converge to $x$, and yet remain ouside $U$, which is impossible.
Therefore, $x_n$ must converge to $x$.\qed\medskip

We'll see later that, under the assumption of metrizability, fundamental sequences, along with their limits, determine all converging sequences in $\bar X$, and thus determine the
compactification (Proposition \ref{C-uniq}).\medskip

We'll now state condiitons that identify the compactification $\tilde X$ we have
built.\medskip

We will say that a certain compactification of $X$ {\em respects
intersections of Weyl chambers} if, given two Weyl chambers based at the point $o$,
$kW$ and $rW$, $k,r\in K$, the intersection of the compactifications of $kW$
and $rW$ is the compactification of the intersection, as in Proposition \ref{props_tildew}.\medskip

The following is a generalization of Lemma 3.18 in \cite{GJT}.

\begin{prop}
\label{C-uniq}
Let $X$ be a locally compact topological space, and take $(i_1,K_1)$, $(i_2,K_2)$ two
metrizable compactifications of $K$. Suppose that ${\cal C}$ is a family of
fundamental sequences for both compactifications (with the possibility that two fundamental sequences may converge to the same limit in one compactification, and to different ones in the other).\m

1. If, for every sequence $(x_n),(y_n)\in \cal C$, $\lim i_1(x_n)=\lim i_1(y_n)$
implies $\lim i_2(x_n)=\lim i_2(y_n)$ then $K_1$ refines $K_2$.\m 

2. If, for every sequence $(x_n),(y_n)\in \cal C$, $\lim i_1(x_n)=\lim i_1(y_n)$ if and
only if $\lim i_2(x_n)=\lim i_2(y_n)$ then $K_1$ and $K_2$ are homeomorphic.
\end{prop}
\proof 1. We will build a continuous map $\phi$ from $K_1$ to $K_2$, which will be a
homeomorphism in the second case.

For $x\in X$, take $\phi(i_1(x)):=i_2(x)$. Now for $x'=\lim x_n$, $(x_n)\in X$, $x'\in \partial
K_1:=K_1\setminus i_1(X)$, let $\phi(x')$ be the common limit in $K_2$ of all sequences in
${\cal C}$ that are subsequences of $(x_n)$ (which belong all to the same class, we can just
take one, and find the limit from that one). The map is clearly onto, and continuous on
$i_1(X)$, we'll now prove continuity at the points
$x'\in\partial K_1$. Notice first that given any point in $x' \in \partial K_1$, there exists a
sequence in
$\cal C$ that converges to it (taking it to be a subsequence of a sequence in $X$ converging to
it, if necessary). 

Let $x' \in \partial K_1$  and suppose there exists a sequence $(x_n)$ of elements of $X$,
with $i_1(x_n)\to x'$. Then we must have $\phi(i_1(x_n))=i_2(x_n)\to \phi(x')$, otherwise, it
would have a subsequence $(y_n)$ not converging to $\phi(x)$. This cannot be,
since any subsequence of $(y_n)$ pertaining to $\cal C$ would converge in $K_2$ to $\phi(x')$,
by definition of $\phi(x')$. 

Now suppose that the sequence $(x_n)$ converging to $x'$ has elements in $\partial K_1$. For
each element $x_n\in \partial K_1$ take $y_n$ to be an element such that both
$d_1(i_1(y_n),x')$, $d_2(i_2(y_n),x')<1/n$, where $d_1$ and $d_2$ are distances in $K_1$ and
$K_2$ respectively; if $x_n\in X$, take $y_n:=x_n$. We have thus built a sequence $(y_n)$ of
elements of $X$ such that $i_1(y_n)\to x'$ and $\lim i_2(x_n)=\lim i_2(y_n)$, if the second
one exists. By the first part of the proof it does exist, and it is equal to $\phi(x')$, thus
$\lim \phi(i_1(x_n))=\lim i_2(x_n)=\lim i_2(y_n)=\phi(x')$.

In case 2, the map is bijective, and the continuity of $\phi^{-1}$ comes from symmetry of
roles of $K_1$ and $K_2$.\qed

2. Any of the two compactifications coincides with the metric
completion of $X$, with respect to the respective metric. This completion is
completely determined by Cauchy sequences in $X$, in either case, and these are
exactly the sequences in $X$ which converge in the metric completion, which is the
compactification. Now, as we have seen (proposition \ref{FundConv}), fundamental
sequences determine the sequences in $X$ that converge in the
compactification. 

Alternatively, the result is also a consequence of theorem 22, p.\ 151, of \cite{Ke}.\qed\medskip

So, briefly put, point 2 in the previous proposition states that, if we have a
metrizable compactification of $X$ admitting a certain class of fundamental sequences,
with a convergence rule, then this is enough to identify the compactification.

\begin{theo}
Suppose that we have a certain metrizable compactification of $W$. Then there is,
up to homeomorphism, at most one compactification of $X$ satisfying the
following properties:
\begin{enumerate}
\item It is metrizable.
\item It is a $K$-compactification.
\item When restricted to $W$ it coincides with the one we have.
\item It respects intersections of Weyl chambers.
\end{enumerate}
Moreover, this compactification is a refinement of any other compactification satisfying
conditions 1--3.\medskip

This compactification exists if the compactification of $W$ is facially stratified. 
\end{theo}
\proof As to existence, the compactification ${\tilde X}$ that we have
constructed before has all the required properties, as we noted in Theorem \ref{compact_metrizable} and Proposition
\ref{props_tildew}.

Now, to check uniqueness, we'll use fundamental sequences. Take two
compactifications $(i_1,X)$, satisfying 1--4, and $(i_2, Y)$, satisfying 1--3. Take the set of
fundamental sequences as in the example:
$$\{k_nx_n.o:k_n {\rm\ and\ } x_n{\rm\ converge} \}.$$
Now we have to prove that equality
of limit in $X$ implies equality of limit in $Y$. We'll just check sequences that converge to
points on the boundary, since for the others, the result is clear.

Suppose then that $k_nx_n$  and $r_ny_n$
are two fundamental sequences with the same limit in $X$, we want them to have the same
limit in $Y$. If $k_n\to k$, $r_n \to r$, then
$\lim kx_n=\lim k_nx_n$, and $\lim ry_n=\lim r_ny_n$, in both
$X$ and $Y$, by continuity of the action of $K$. 

The common limit point in $X$ is thus in the compactification of $kW$ and
$rW$. By condition 3, there must exist a sequence in $(z_n)$ in $kW\cap rW$, converging to
the same point. Now, by $K$-equivariance,
$$\lim i_1(z_n)=\lim i_1(kx_n)\Rightarrow \lim i_1(k^{-1}z_n)=\lim i_1(x_n),$$
and the last limits are in $W$. Since the compactification of
$W$ coincides in both $i_1$ and $i_2$, then $\lim i_2(k^{-1}z_n)=\lim i_2(x_n$, and
$\lim i_2(z_n)=\lim i_2(kx_n)$. Similarly, $\lim i_2(z_n)=\lim i_2(ry_n)$, and
thus the limits are the same in $Y$. This proves that $X$ refines $Y$, by proposition
\ref{C-uniq}. 

Now if $Y$ satisfies also condition 4, we can repeat the argument with $X$ and $Y$
interchanged. We thus get that limits of fundamental sequences coincide on $X$ and $Y$,
and this proves that the compactifications are homeomorphic, again by proposition
\ref{C-uniq}.\qed\medskip

{\bf Example}. Take the symmetric space ${\mathbb H}^2\cong \sl(2,\R)/\so(2)$. Here the
dominant Weyl chamber is not more than a half-geodesic starting from $o$. We can
compactify it by joining a point to it, and there will be at least two $K$-equivariant
compactifications that restricted to $W$ will be this one: the one-point
compactification (adding a point $\infty$ to ${\mathbb H}^2$) and the visual
compactification. However, the one-point compactification does not respect
intersections of Weyl chambers, since the intersection of two Weyl chambers is
$\{o\}$, and the intersection of their compactifications is $\{o, \infty\}$.

As we see, the visual compactification, which respects intersections of Weyl
chambers, refines the one-point compactification.

\section{Addendum: Generalized Busemann Compactifications}

We now present another way of building compactifications of symmetric spaces, which generalizes Busemann compactifications. We will not explore this concept as much as the previous one, but limit ourselves to proving that it does indeed produce a compactification.\medskip

We start with a function $\delta:X\times X\to C$, where $C\subseteq \R_+^n$ is a convex cone, and we'll assume this function has the following properties:
\begin{enumerate}
\item Its norm should be strictly increasing with distance, i.e.\ if $d(z,x)>d(y,x)$, then $||\delta(z,x)|| > ||\delta(y,x)||$, with $\delta(x,x)=0$. 
\item A Lipschitz condition: for some $s>0$,
$$||\delta(x,y)||\leq s d(x,y),$$
\item A triangle inequality: for some $k$, $||\delta(x,y)-\delta(x,z)||\leq k d(x,z)$.
\end{enumerate}

We now prove that, under these conditions, this function (which we can call a {\em kernel}) can be used to define a compactification of $X$ in the same way the distance function is used to define the Busemann compactification. 

To this end, fix a point $o\in X=G/K$, and, for a given $x\in X$, define $b_{x}:X\to C$ as 
$$b_{x}(y):=\delta(x,y)-\delta(x,o).$$ 
Taking in $\cont (X, C)$ the topology of uniform convergence on compacts, we now show that the map $\phi:x\mapsto b_{x}$ is an embedding of $X$ in $\cont (X,C)$, using then the Ascoli-Arzel\`a theorem to prove its image is compact. We'll follow \cite{Bal}, but with a different notation.\medskip

We first prove that, for a given $x\in X$, $b_{x}$ is Lipschitz. Given $z,z'\in X$, we have 
\begin{eqnarray*}
||b_{x}(z)-b_{x}(z')|| & = & || \delta(x,z)-\delta(x,o)-\delta(x,z')+\delta(x,o)||\\
& = & || \delta(x,z)-\delta(x,z')||\  \leq \  k d(z,z')
\end{eqnarray*}

Therefore, the function $\phi$ maps $X$ to $\cont (X,C)$. Now, to see it is one to one, take $x\neq x'\in X$. Because $C\in \R_+^n$, we must have that either $d(x,o)-d(x',o)\not\in C$ or $d(x',o)-d(x,o)\not\in C$, assume the first case holds. Then, we have
\begin{eqnarray*}
b_x(x')-b_{x'}(x') & = & \delta(x,x')-\delta(x,o)-\delta(x',x')+\delta(x',o)\\
& = & \delta(x,x')-(\delta(x,o)-\delta(x',o))\neq 0,
\end{eqnarray*}
which proves the map is one to one. \vs

To check that it is an embedding, suppose that $b_{x_n}\to b_x$, and $x_n\not\to x$. If the sequence $x_n$ remains bounded, then it must have a converging subsequence, and by what we already proved, this subsequence has to converge to $x$. Since this has to be true of any converging subsequence, we have the result in this case. 

We now consider the case where $x_n$ is not bounded. In this case, consider the closed ball of radius 1 around $x$, $B$,  which is a compact set. We must have that $||b_{x_n}||\to ||b_x||$ inside the ball. Consider, for each $n$, the geodesic going from $x_n$ to $x$. The function 
$$||b_{x_n}||=||\delta(x_n,\cdot)-\delta (x_n,o)||$$ 
has to be increasing, along this geodesic, as we move from $x_n$ to $x$, because of condition 1, but $||b_x||$ has a minimum at $x$, $b_x(x)=0$. Denoting by $\partial B=\{y\in X :  d(y,x)=1\}$, let $m=\min_{y\in \partial B} ||b_{x}(y)||$. By condition 1, we have $m > 0$. Take $0<\epsilon<m/2$ and the set
$$V_\epsilon : = \{f\in C:\forall y\in b\ ||f(y)||-||b_x(y)|| |<\epsilon \}.$$

For functions in $V_\epsilon$ and $y\in B$, 
$$|\, ||f(y)||-||b_x(y)||  \, | \leq ||f(y)||-||b_x(y)||<\epsilon.$$
So given any function $f\in V_\epsilon$ and $y\in \partial B$, we must have $f(y)>\epsilon$, and $f(x)<\epsilon$. If we had any $b_{x_n}$ inside this neighborhood, its restriction to the geodesic from $x_n$ to $x$ could not be an increasing function, because of the previous considerations.

Hence, we can't have the uniform convergence in this ball. \medskip

Under these conditions, the Ascoli-Arzel\`a theorem assures that the closure of the set $\{b_x:x\in G/H\}$ is a compact set, yielding therefore a compactification of the symmetric space. We thus have proved the following result.

\begin{theo}
Let $\delta$ be a function satisfying conditons 1.-3.\ above, and find a point $o\in X$. Consider the map 
\apl{X}{\cont (X,C)}{x}{\delta(x,\cdot)-\delta(x,o)}
Then the closure of the image of $X$ in $\cont (X,C)$ is a compactification of $X$.
\end{theo}

We now present some functions that we conjecture will yield the known compactifications that we have mentioned. 

For $x=gK$ and $y=hK$, define $r(x,y)$ as the {\em generalized radius of $y$ from $x$}, which can be defined as the element $H\in \overline{\a^+}$ such that $g^{-1}hK= k e^H K$, for some $k\in K$. It's easy to check that this is well defined, and that it coincides with the usual generalized radius of $y$ if we choose $x$ as a reference point in the symmetric space instead of $o$. 

Consider that the set of simple roots is ordered, and for $I\subseteq \Delta$, consider, for $H\in \overline{\a^+}$, $(\alpha(H):\alpha\in I)$ as a well defined element of $(\R_0^+)^{|I|}$. We denote this element by $\alpha_I(H)$. 

We now present the functions that we claim yield the compactifications we studied. 
\begin{itemize}
\item The function $\delta(x,y)=(\alpha(r(x,y)): \alpha\in \Delta)$ yields the maximal Furstenberg/maximal Satake compactification.
\item The function $\delta(x,y)=((||\alpha_I(r(x,y))||): |I|=1 \text{ or } 2)$ yields the Martin compactification.
\item The function $\delta(x,y)=((||\alpha_I(r(x,y))||): I\subseteq \Delta)$ yields the Karpelevi\v c compactification.

\end{itemize}

\end{document}